\def\R{{\rm I}\!{\rm R}}

\def\carre{\hbox{\vrule \vbox to 7pt{\hrule width 6pt \vfill \hrule}\vrule }}

\font\larger = cmr10 at 14pt

\font\small = cmr10 at 8pt

\parindent = 0cm

\centerline { \larger
On bounded pseudodifferential operators  in a high-dimensional setting}

\medskip

\centerline { L. Amour, L. Jager,  J. Nourrigat}

\medskip

\centerline { Universit\'e de Reims}

\medskip

\centerline {\it Dedicated to the memory of Bernard Lascar}

\vskip 1cm

\centerline{\small \vbox{\hsize 16 cm{\noindent{\bf Abstract}}\vskip 0.1 cm
{ This work is concerned with extending the results of Calder\' on and
 Vaillancourt proving  the boundedness of
  Weyl pseudodifferential operators
$ Op_h^{Weyl} (F)$
 in
$L^2(\R^n)$. We state conditions under which the norm of such operators
has an upper bound independent of  $n$.
To this aim, we apply a decomposition of the identity to the symbol
$F$, thus obtaining a sum of operators of a hybrid type, each of them
 behaving as a
  Weyl
 operator with respect to some of the variables and as an anti-Wick operator
 with respect to the other ones.
Then  we establish upper bounds for  these auxiliary operators,
using suitably adapted classical methods like coherent states. }
}}

\bigskip

 {\sevenrm
2010  Mathematical Subject Classification   35S05

  }
\bigskip

{\bf 1. Introduction.}

\bigskip

Since the work of  Calder\'on and Vaillancourt [C-V], it is well known that,
if a function $F$, defined on  $\R^{2n}$, is smooth and has bounded
 derivatives,
 it is possible to associate with it a   pseudodifferential
 operator, depending on a parameter $h>0$,  which is bounded on $L^2 (\R^n)$ (see also [HO], [LER], [R], [U2]).
 This operator is formally defined by:
$$ (Op_h^{Weyl} (F) f) (x) = (2 \pi h)^{-n} \int _{\R^{2n}} e^{ {i \over h}(x-y)\cdot \xi }
 F \left ( { x+y \over 2} , \xi \right )f(y) dyd\xi \hskip 2cm  x\in \R^n \leqno (1.1)$$
for $f$ belonging to $L^2(\R^n)$. (When $h=1$ the subscript $h$ will
be omitted). Moreover, its norm is bounded above by
$$ \Vert Op_h^{Weyl} (F) \Vert _{{\cal L}(L^2(\R^n)) } \leq C
\sum _{|\alpha +\beta|\leq N} \Vert
\partial_{x}^{\alpha} \partial_{\xi} ^{\beta}F \Vert _{L^{\infty}
(\R^{2n})} \leqno (1.2)$$
where $N$ and $C$ depend on the  dimension $n$.

The aim of this work is to prove that, under certain conditions,
the constants appearing in the upper bound do not depend on the dimension.
 The set of derivation  multi-indices which are used  depends
 on the dimension in a way that will be precisely stated.

We shall thus be able to give examples where the dimension goes
to infinity and the norm, nevertheless, remains bounded.
\bigskip
In a later work we shall study pseudodifferential operators
where the configuration space  $\R^n$  will be replaced by
an infinite dimensional Hilbert space, by a method differing from
Bernard Lascar's (see [LA1]- [LA10]).
These results have been announced in a preprint [A-J-N] in September 2012.

\bigskip

We first recall  an example in which the constant appearing in
the upper bound on the norm does not depend on the dimension. This is the case
 when
the function $F$ is the Fourier transform of a function $G$ belonging to
$L^1(\R^n)$ :
$$ F(x , \xi) = (2\pi h)^{-2n} \int _{\R^{2n}} e^{ {i\over h}
(a\cdot x + b\cdot \xi) } G(a , b) dadb .$$
Since the Weyl operator associated with the function
$$(x , \xi) \rightarrow  E_{a , b, h}(x , \xi) = e^{{i\over h} (a \cdot x + b\cdot \xi )}$$
is the operator $W_{a , b , h}$ defined by
$$ (Op_h^{Weyl} (E_{a , b, h}) f ) (u) = (W_{a , b , h}f) (u) =
 e^{ {i\over h} a \cdot u + { i\over
2h} a \cdot b } f(u + b), $$
 the equality  (1.1) may be rewritten in the form
$$ Op_h^{Weyl} (F) =(2\pi h)^{-2n} \int _{\R^{2n}} G(a , b)W_{a , b , h} dadb . $$
Since $W_{a , b , h}$ is unitary,
 $$ \Vert Op_h^{Weyl} (F) \Vert \leq (2\pi h)^{-2n} \int _{\R^{2n}} |G(a ,
 b)|\ dadb .$$
Situations of this kind have been considered by B. Lascar ([LA1]-[LA10]) in
 an infinite dimensional setting, but the $L^2$ boundedness was not
 the main motivation of these works.

\bigskip

Our approach is different, in that we aim at extending the bound (1.2).
 Let us specify the
set of multi-indices which will be used.
 Cordes [C], Coifman Meyer [C-M], Hwang [HW] noticed that one does not need
all the multi-indices to state (1.2) but only  the
  $(\alpha , \beta )$ satisfying
 $0\leq \alpha_j \leq 1$ and $0\leq \beta_j \leq 1$ for each $j$.
In this paper we shall use the multi-indices   $(\alpha ,
\beta )$ such that $0\leq \alpha_j \leq 2$ and $0\leq \beta_j \leq 2$
for each $j$.
We now can state the hypotheses on the function $F$.

\bigskip

Let  $(\rho_j )_{1\leq j \leq n}$
 and  $(\delta_j )_{1\leq j \leq n}$ be two sequences satisfying
 $\rho_j \geq 0$ and $\delta_j \geq
0$ for every $j\leq n$, let $M$ be a
nonnegative real number.
Suppose that

\bigskip

(H) for every multi-index $(\alpha , \beta)$ such that  $0\leq
\alpha_j \leq 2$ and $0 \leq \beta_j \leq 2$ for every $j\leq n$,
the partial derivative $\partial_x^{\alpha}
\partial_{\xi}^{\beta} F$ exists, is continuous, bounded and satisfies

$$ |\partial_x^{\alpha} \partial_{\xi}^{\beta} F(x , \xi)| \leq M
 \prod _{j=1}^n \rho_j ^{\alpha _j} \delta _j^{\beta_j}. \leqno (1.3)$$
If  $\rho_j = 0$ and  $\alpha_j= 0$, we set that
 $\rho_j^{\alpha_j}=1$.

\bigskip

Our main result is the following Theorem.

\bigskip

{\bf Theorem 1.1.} {\it If a function $F$ defined on $\R^{2n}$
satisfies  hypothesis (H), then the operator $Op_h^{Weyl}(F)$, defined
 formally by  (1.1), is bounded in
$L^2(\R^n)$ and satisfies
$$ \Vert Op_h^{Weyl} (F) \Vert _{{\cal L} (L^2(\R^n))} \leq M \
\prod _{j=1}^n (1 + 81 \pi  h \rho_j \delta_j  ) \leqno (1.4)$$
 if $0 < h \rho_j \delta_j \leq 1$ for every $j\leq n$.
}

\bigskip

{\it Example 1.2.} Let $V\geq 0$ be a real-valued bounded function
in $C^{\infty}(\R)$, whose derivatives are all bounded. For all
integer $n\geq 1$, set
$$ H_n (x , \xi) = \sum  _{j\leq n} \xi_j^2 + \sum _{j \leq n,
k\leq n \atop |j-k| = 1} g_j g_k V(x_j - x_k) $$
where $(g_j)$ is a sequence of positive numbers such that, for some
$C_0>0$, we have $g_j \leq C_0 g_k$ if $|j-k|\leq 1$. Set:
$$ P_n (x , \xi) = e^{- H_n (x , \xi)} \leqno (1.5)$$
 We shall see that Hypothesis (H) is
satisfied, with
$$M = 1, \hskip 2cm \delta_j = C_1 \hskip 2cm   \rho_j = C_1 \lambda_j
\leqno (1.6) $$
$$ \lambda_j =  \max _{1\leq \nu \leq 4} ( g_j^2 \Vert V^{(\nu)}\Vert _{L^{\infty}} )^{1/\nu}
 \leqno (1.7)$$
where $C_1$ is a real constant, to be determined,  depending only on
$C_0$.   Let us set:
$$ W(x) = -\sum _{j \leq n,
k\leq n \atop |j-k|=  1} g_j g_k V(x_j - x_k), $$
and
$$A_n = \{ (j , k), |j-k|=1, \ \ 1\leq j, k \leq n\} . $$
We shall estimate $ \partial ^{\alpha} e^W $ when $\alpha_j \leq 2$
for every $j$. For each function $f\in C^4(\R)$ with bounded
derivatives, set:
$$ M(f) = \max _{1\leq \nu \leq 4} (\Vert f^{(\nu)} \Vert
_{L^{\infty}} )^{1/\nu}.$$
We notice that:
$$ e^{-f} \partial ^{\nu} e^f = (\partial + f')^{\nu} \cdot 1.$$
From this and a simple computation follows that:
$$ |e^{-f} \partial ^{\nu} e^f| \leq (2 M(f))^{\nu}. \leqno (1.8) $$
In order to apply this inequality we divide  $W$ into two parts.
Set:
$$ U_j = U_j (x_1 , ... , x_n) = - g_j g_{j+1} U ( x_j - x_{j+1}), \
\ \ \ U(y) = V(y)+ V(-y)$$
and
$$ W_e = \sum _{j \ even } U_j, \hskip 2cm W_o= \sum _{j \ odd}
U_j. $$
Then $W = W_e + W_o$ and we notice that the variable $x_{\nu}$
occurs only once in the $W_e$ and $W_o$. Also, since
$$ \partial^{\alpha} e^{W_e} = \prod_{ j\ even, \ j <n} (\partial
_j^{\alpha_j}\partial _{j+1}^{\alpha_{j+1}}e^{U_j} )$$
we have the estimate
$$|\partial^{\alpha} e^{W_e}|\leq e^{W_e} \prod _{j \ even, \ j<
n}T_j \leqno (1.9)$$
where $T_j$ is the $L^{\infty}$ norm of the function:
$$ e^{f_j} \partial ^{\alpha_j+\alpha_{j+1}} e^{-f_j} , \hskip 1cm
f_j = g_j g_{j+1} U(x).$$
Let $\lambda_j$ be defined by (1.7), and set:
$$ M_j =  \max _{1\leq \nu \leq 4}( \Vert g_j g_{j+1}
U^{(\nu)} \Vert _{L^{\infty }})^{1/\nu}.$$
Then
$$ M_j \leq 2C_0 \lambda_j, \hskip 2cm M_j\leq 2 C_0 \lambda
_{j+1}\leqno (1.10)$$
where $C_0$ is defined so that $g_j g_{j+1}\leq C_0 \min (g_j^2,
g_{j+1}^2)$. It follows from (1.8) and (1.10) that:
$$ T_j \leq (2 M_j)^{ \alpha_j+\alpha_{j+1}} \leq ( 4C_0
\lambda_j)^{\alpha_j} (4C_0 \lambda_{j+1}) ^{\alpha_{j+1}}.$$
Then (1.9) gives:
$$ |\partial^{\alpha} e^{W_e} |\leq e^{W_e} \prod _{ j \ even, j<n}
( 4C_0 \lambda_j)^{\alpha_j} (4C_0 \lambda_{j+1}) ^{\alpha_{j+1}}. $$
In a similar way one gets the estimate:
$$ |\partial^{\alpha} e^{W_o} |\leq e^{W_o} \prod _{ j \ odd, j<n}
( 4C_0 \lambda_j)^{\alpha_j} (4C_0 \lambda_{j+1}) ^{\alpha_{j+1}}. $$
Then we write:
$$ \partial^{\alpha} e^W =\partial^{\alpha} e^{W_e}e^{W_o}= \sum '
\pmatrix { \alpha \cr  \beta \cr } (\partial ^{\alpha - \beta}
e^{W_e})(\partial ^{ \beta} e^{W_o})   $$
where the prime indicates that one only takes the summation over
terms with $\beta _1 = \alpha_1$ and $\beta _n = \alpha_n$ if $n$ is
even, and with $\beta _1 = \alpha_1$ and $\beta _n = 0$ if $n$ is
odd. We get the estimate:
$$ |\partial^{\alpha} e^W | \leq e^W  \sum '
\pmatrix { \alpha \cr  \beta \cr }  \Big ( \prod _{j\ even, \ j<n}
(4C_0 \lambda_j) ^{\alpha_j - \beta_j} (4C_0 \lambda_{j+1})
^{\alpha_{j+1} - \beta_{j+1}}\Big ) \cdot $$
$$ \cdot \Big ( \prod _{j\ odd , \ j<n}
(4C_0 \lambda_j) ^{ \beta_j} (4C_0 \lambda_{j+1}) ^{
\beta_{j+1}}\Big ) $$
$$ =e^W  \sum '
\pmatrix { \alpha \cr  \beta \cr }  \prod _{j=1}^n (4C_0
\lambda_j)^{\alpha_j} \leq 2^{|\alpha|} e^W\prod _{j=1}^n (4C_0
\lambda_j)^{\alpha_j}. $$
We have proved that:
$$ |\partial ^{\alpha}e^W|\leq e^W \prod _{j=1}^n (8 C_0
\lambda_j)^{\alpha_j}. $$
Therefore, hypothesis (H) is satisfied with the choice (1.6). We may
apply  Theorem 1.1 if $h \rho_j \delta_j \leq 1$ for all $j$. It
follows that $\Vert Op_h^{Weyl}(P_n) \Vert _{{\cal L}(L^2
(\R^{n}))}$  is bounded independently of $n$ if the sum $\sum
_{j\geq 1} g_j^{1/2}$ converges, and if $h$ is small enough.
 If $g_j = 1$, the norm is not bounded, but estimated, with some
 constant $C>0$, independent of the dimension, by:
$$ \Vert Op_h^{Weyl}(P_n) \Vert _{{\cal L}(L^2 (\R^{n}))} \leq
e^{ C h n }  \leqno (1.11). $$

\bigskip

{\it Example 1.3.}
The mean-field approximation  uses hamiltonians
 of the form
$$ H_n (x , \xi) = \sum  _{j\leq n} \xi_j^2 +{1\over n}  \sum _{j \leq n,
k\leq n }  V(x_j - x_k) $$
where  $V$ is as in Example 1.2. Let $P_n$ be the function defined
as in (1.5). Then hypothesis (H) is satisfied with $M=1$ and $\rho_j
= \delta _j = C_1 $, where $C_1$ does not depend on $n$. In this
case, Theorem 1.1 shows that, provided $C_1$ is small enough, we
have also (1.11) for some constant $C$ which is independent of $n$.

\bigskip

We express our thanks to the referee for his helpful suggestions, which
allowed us in particular to gain on the number of derivatives  and to
simplify the proofs.

\bigskip
{\bf 2. Hybrid Weyl anti-Wick  quantization.}

\bigskip

In order to prove Theorem 1.1, we may as well assume that
$\rho_j = \delta _j $ for every  $j\leq n$ and that  $h=1$.
Indeed, if a function  $F$ satisfies  hypothesis (H) with two
sequences $(\rho_j)$ and $(\delta_j)$ of positive real numbers,
then the function $\widetilde F$ defined by
$$ \widetilde F (x , \xi)= F \left ( x_1 \sqrt h \lambda_1 , ... ,
 x_n \sqrt h \lambda_n , {\xi_1 \sqrt h \over \lambda _1 } , ... ,
 {\xi_n \sqrt h \over \lambda _n } \right ) \hskip 2cm \lambda _j = \sqrt { \delta _j \over \rho_j}$$
satisfies (H) with  $\rho_j$ and $\delta_j$ replaced by
 $ \varepsilon_j = \sqrt {h \rho_j
\delta _j}$.
If Theorem 1.1 is valid for  $Op_1^{Weyl}$,
then we  get that, if  $\varepsilon_j ^2 \leq  1$:
$$ \Vert Op_1^{Weyl} (\widetilde F) \Vert _{{\cal L} (L^2(\R^n))}
\leq M \ \prod _{j=1}^n (1 + 81 \pi  \varepsilon _j^2 ) .$$
Since $Op_h^{Weyl} (F) = T^{-1} Op_1^{Weyl } (\widetilde F) T$,
where $T$ is a unitary operator acting in $L^2(\R^n)$, Theorem 1.1
for  $Op_h^{Weyl} (F)$ holds true. This  follows by continuity  if
some of the  $\rho_j$ or  $\delta_j$ are equal to $0$.

\bigskip

Consequently,  we shall assume from now on that
 $\rho_j =
\delta_j \leq 1$ for every  $j\leq n$ and that $h=1$.

\bigskip

In order to prove Theorem 1.1 we shall split the operator into
a sum of operators which will behave as Weyl operators with respect to
a first subset of the variables (meaning the
operators will be defined by a
 formula analogous to (1.1)  in which only these variables appear)
and as anti-Wick operators with respect to the other variables.

\bigskip

We first need to recall the anti-Wick quantization. The definition
uses the coherent states, which is the  family of functions $
\Psi_{X } $ indexed by   $X = (x , \xi)\in \R^{2n}$,
  defined by
$$ \Psi_{X } (u) =  \pi ^{ -n/4} e^{-{| u-x|^2 \over 2}} e^{i  u .\xi -
{i \over 2} x. \xi} \hskip 2cm X = (x , \xi) \in \R^{2n} \hskip 1cm
u\in \R^n. \leqno (2.1) $$
Recall that
$$  < f , g> = (2\pi )^{-n} \int _{\R^{2n}} < f , \Psi_{X }>
\ < \Psi _{X } , g> dX . \leqno (2.2) $$
If  $F$ is a function in $L^{\infty}(\R^{2n})$,
one can associate with it an  (anti-Wick
)  operator
 $ Op^{AW}(F) $ such that, for all  $f$ and  $g$ in  $L^2(\R^n)$:
$$ < Op^{AW}(F) f , g> = (2\pi )^{-n} \int _{\R^{2n}}F(X)  < f , \Psi_{X }>
\  < \Psi _{X }, g> \  dX  . \leqno (2.3) $$
We then have
$$ \Vert Op^{AW} (F)\Vert _{{\cal L} (L^2(\R^n))}
 \leq  \Vert F\Vert _{L^{\infty}(\R^{2n})} . \leqno (2.4)$$
The relationship between Weyl and anti-Wick quantizations is given, for every
 $F$ in $L^{\infty}(\R^{2n})$, by :
 $$ Op ^{AW }(F) = Op^{Weyl} \left ( e^{{1\over 4} \Delta
} F \right ) \leqno (2.5) $$
where
$$\Delta  = \sum _{j\leq n} \Delta_j \hskip 2cm \Delta_j =
{\partial^2\over \partial x_j^2} + {\partial^2\over \partial
\xi_j^2} . \leqno (2.6) $$
This fact is classical  (see Folland [F]).
 One has an identity decomposition in
 $L^{\infty}(\R^{2n})$:
$$I = \sum _{E \subseteq \{1 , ..., n\} } T(E)e^{{1\over 4} \Delta
_{E^c}} \hskip 2cm T(E) = \prod _{j\in E} (I -e^{{1\over 4} \Delta
_j}) \leqno (2.7) $$
$$\Delta _{E^c} =\sum _{j\in E^c} \Delta_j .\leqno (2.8) $$

\bigskip

For every subset  $E \subseteq \{ 1 , ... , n\}$ and every symbol
$F$, we define an operator $ Op ^{hyb , E}(F)$ by :
$$ Op ^{hyb , E}(F) = Op^{Weyl} \left ( e^{{1 \over 4} \Delta
_{E^c}} F \right ) \leqno (2.9) $$
This operator behaves as a Weyl operator with respect to the
variables
 $x_j$ $(j\in E)$ and as an anti-Wick operator with respect to the
 variables $x_j$ $(j\in E^c)$. If $E = \emptyset$, it is the
 anti-Wick operator and conversely if
 $E = \{ 1 , ... , n\}$, it is the Weyl operator.

\bigskip

One derives a decomposition of the Weyl operator $Op^{Weyl } (F)$:
$$ Op^{weyl } (F) =\sum _{E \subseteq \{1 , ..., n\} } Op^{hyb,
E} (T(E)F) .\leqno (2.10) $$

\bigskip

We shall now prove an upper bound on the norm of a hybrid operator
$Op^{hyb, E} (G)$, where the function  $G$ is bounded on $\R^{2n}$.
The only derivatives of $G$ which will play a role are the
derivatives with respect to $x_j$ or $\xi_j$ with $j\in E$.
 For every integer $m$ we introduce the set of multi-indices
$$I_m(E) = \{ (\alpha , \beta ) , \ \ \ \alpha _j \leq m , \ \ \ \
\beta_j \leq m, \ \ (1\leq j \leq n) \ \ \  \ \ \ \    \alpha_j =
\beta _j = 0 \ \ \ \ \ {\rm if} \ \ j\notin E \}.  \leqno (2.11)$$

\bigskip

We shall prove the following Lemma in Section 3, by adapting classical methods
(Unterberger [U2]).

\bigskip

{\bf Lemma 2.1.} {\it If $F$ satisfies hypothesis (H) and if $E
\not= \emptyset$, then
$$\Vert Op^{hyb, E} (F) \Vert_{{\cal L}(L^2(\R^n))}
 \leq  \left ( {9 \pi \over  2} \right )^{|E|} \sum _{(\alpha ,
\beta)\in I_2(E)}\Vert
\partial_x^{\alpha } \partial _{\xi}^{\beta} F\Vert
_{L^{\infty} (\R^{2n})} \leqno (2.12) $$

}

\bigskip

We shall establish the following Lemma in Section 4.

\bigskip
 {\bf Lemma 2.2.} {\it If $F$  satisfies hypothesis  (H) with
$\rho_j = \delta_j \leq 1 $ for every $j \leq n$, and if $E \not=
\emptyset$, the function $T(E)F$   satisfies
$$  \sum _{(\alpha , \beta)\in I_2(E)}
\Vert
\partial_x^{\alpha } \partial _{\xi}^{\beta} T(E) F\Vert
_{L^{\infty} (\R^{2n})} \leq M 18 ^{|E|}\prod _{j\in E}\rho_j^2
\leqno (2.13) $$
where  $T(E)$ is defined in (2.7). }

\bigskip

{\it End of the proof of Theorem 1.1.} We assume that $h=1$ and that
 $\rho_j = \delta_j \leq 1$ for every $j\leq n$. According to
(2.10), we have
$$ \Vert Op^{Weyl} (F) \Vert _{{\cal L} (L^2(\R^n))}
 \leq \sum _{E \subseteq \{1 , ..., n\} } \Vert Op^{hyb, E}
(T(E) F) \Vert_{{\cal L} (L^2(\R^n))}.
$$
By Lemma 2.1:
$$ \Vert Op^{Weyl} (F) \Vert _{{\cal L} (L^2(\R^n))}
 \leq \sum _{E \subseteq \{1 , ..., n\} } \left (  {9 \pi \over  2} \right )^{|E|} \sum _{(\alpha ,
\beta)\in I_2(E)}  \Vert
\partial_x^{\alpha } \partial _{\xi}^{\beta} T(E)F\Vert
_{L^{\infty} (\R^{2n})}.
$$
With the same hypotheses, Lemma 2.2 shows that:
$$ \Vert Op^{Weyl} (F) \Vert _{{\cal L} (L^2(\R^n))} \leq M \
\sum _{E \subseteq \{1 , ..., n\} } \left (  {9  \pi  \over 2}
\times 18   \right )^{|E|} \prod _{j\in E} \rho_j^2. $$
It follows easily that
$$ \Vert Op^{Weyl} (F) \Vert _{{\cal L} (L^2(\R^n))} \leq M \prod _{j\leq n} ( 1 + 81 \pi   \rho_j^2). $$
The theorem is proved in the case when  $h= 1$ and $\rho_j =
\delta_j \leq 1$ for all $j\leq n$. In the general case, the
announced result follows as we saw.

\bigskip

{\bf 3. Proof of Lemma 2.1.}

\bigskip

We shall use the results of Unterberger [U1], [U2] concerning the
upper bound of
 $< A \Psi_{X } , \Psi_{Y }>$, where  $A$ is a pseudodifferential operator
and the
  $\Psi_{X }$ are the coherent states defined by (2.1).
We first recall the integral expression of this scalar product and give an analogous
statement for hybrid operators.

\bigskip

{\bf Proposition 3.1.} {\it Let  $F$ be a function defined on  $\R^{2n}$ and
 satisfying
 hypothesis  (H). Then we get, for every  $X$ and $Y$ in
$\R^{2n}$:
$$ <Op^{Weyl}(F )
\Psi _{X} , \Psi_{Y} > =  \pi ^{-n}\int _{\R^{2n}}
 F(Z)\   \Phi_{n } (X , Y , Z)  d Z \leqno (3.1)$$
with
$$  \Phi_{ n } (X , Y , Z) = e^{ -|Z -{X+Y\over 2}|^2 - i \sigma
(Z , X-Y) - {i\over 2}\sigma (X , Y) }  \leqno (3.2) $$
where the symplectic form  $\sigma$ is given by  $\sigma (X , Y) =  y
\cdot \xi - x \cdot \eta$ for all  $X = (x, \xi)$ and  $Y = (y , \eta)$ in  $\R^{2n}$,
}

\bigskip

{\it Proof.}
For all functions $f$ and $g$ belonging to the Schwartz space  ${\cal S}(\R^n)$,
one defines the Wigner function
  $H (f , g , Z)$ ($Z \in \R^{2n}$) by :
 $$H (f , g , Z) = \int _{\R^n } e^{-i t\cdot \zeta } f \left (
 z + {t\over 2} \right ) \overline {g \left (
 z - {t\over 2} \right )} dt \hskip 2cm Z = (z , \zeta) \in \R^{2n} \leqno (3.3)$$
 (cf Unterberger [U2], or Lerner [LE], sections 2.1.1 et 2.1.2,
  or   Combescure  Robert [C-R], section 2.2).
The following equality is proved in [U2] or [LE] or [C-R],
 for all  $f$ and $g$
in  ${\cal S}(\R^n)$ and every Borel function  $F$ which is bounded on $\R^{2n}$:
 $$  <  Op^{Weyl}(F)  f ,  g>  =(2\pi ) ^{-n}
  \int _{\R^{2n}} F(Z)H (f , g , Z) dZ .  \leqno (3.4)$$
An explicit computation using the coherent spaces  $\Psi_{X }$  defined by (2.1)
shows that
$$ H (\Psi_{X } , \Psi_{Y } ,
Z) = 2^{n } \Phi_{n } (X , Y , Z) , \leqno (3.5)$$
which implies (3.1). \hfill $\carre$

\bigskip

 Let $n' < n$ and  $n'' = n - n'$. We denote by
  $X = (X' , X'')$ the variable in $\R^{2n}$, with  $X' = (X_1 ,
... , X_{n'})$ and  $X'' = (X_{n' + 1} , ... , X_n)$. Set
$$ \Delta '' = \sum _{ j = n'+1}^n ( \partial_{x_j}^2 +
\partial_{\xi_j}^2) .$$

\bigskip

{\bf Proposition 3.2.} {\it For all $f$  and $g$ in  ${\cal
S}(\R^n)$, we have :
$$ < Op ^{Weyl} ( e^{{1\over 4}\Delta ''} F) f, g> = ... \leqno (3.6) $$
$$ = C(n', n'')  \int _{\R^{6n'+ 2n''}} F
(Z', T'') \Phi_{n' } ( X', Y', Z') < f, \Psi_{X', T'' }> < \Psi_{Y',
T'' } , g> dX'dY'dZ' dT''  $$
where  $ C(n', n'') = 2^{n'} (2\pi )^{-3n' - n''}$,  $\Phi_{n' }$
being the function defined by (3.2) with $n',X'$ instead of $n,X$.
}

\bigskip

{\it Proof.} Proposition 3.1 and the
integral expression for the heat operator $ e^{{1\over 4}\Delta ''}$
give:

$$ < Op ^{Weyl} ( e^{{1\over 4}\Delta ''} F) \Psi_{X } , \Psi_{Y }> =
\pi ^{-n- n''} \int _{\R^{2(n+n'')}} e^{-|Z''-T''|^2 } F( Z', T'')
 \Phi_{n } (X , Y , Z) d ZdT'' \leqno (3.7)$$
An explicit computation yields
$$ \pi ^{-n''} \int _{\R^{2n''}} e^{-|Z''-T''|^2 }
 \Phi _{n'' } (X'' , Y'' , Z'') dZ'' =2^{-n''}
< \Psi _{X'' } , \Psi_{T'' }> _{n''} < \Psi_{T '' } , \Psi_{Y''
}>_{n''}
$$
where  $< \cdot , \cdot >_{n''}$ is the scalar product of
$L^2(\R^{n''})$.  Combining (3.7) with this equality one sees that
$$ < Op ^{Weyl} ( e^{{1\over 4}\Delta ''} F) \Psi_{X } , \Psi_{Y}> = ... \leqno (3.8) $$
$$ =2^{-n''} \pi ^{-n} \int _{\R^{2n}} F
(Z', T'') \Phi_{n' } ( X', Y', Z') < \Psi _{X'' } , \Psi_{T''
}>_{n''} < \Psi_{T'' } , \Psi_{Y'' }>_{n''} dZ' dT'' . $$
For all  $f$ and $g$ in ${\cal S}(\R^n)$, one gets,  applying (2.2) twice :
$$ < Op ^{Weyl} ( e^{{1\over 4}\Delta ''} F) f, g> = ... \leqno (3.9) $$
$$ (2 \pi )^{-2n}\int _{\R^{4n}} < f , \Psi _{X  }> < Op ^{Weyl}
( e^{{1\over 4}\Delta ''} F) \Psi_{X } , \Psi_{Y }> < \Psi_{Y } , g>
dX dY  . $$
One then applies  (3.8) and the following result,
deduced from  (2.2) in  dimension $n''$:
$$ (2 \pi)^{-n''}\int _{\R^{2n'' }}< f ,\Psi _{X  }> < \Psi _{X''
}, \Psi _{T'' }>_{n''}  d X'' = < f ,\Psi _{X' , T''  }> . $$
Formula (3.6) follows from that and from an analogous result
about
$<\Psi_{Y } , g>$.  \hfill $\carre$

\bigskip

 For every $X' = (x' , \xi' )$ in
$\R^{2n'}$, set:
$$ K_{n' } (X') = \prod _{j= 1} ^{n'} \left ( 1 +  x_j^2  \right )
 \left ( 1 + \xi_j^2  \right ) .\leqno (3.10)$$

\bigskip

{\bf Lemma 3.3.} {\it For every function $G$ satisfying
hypothesis  (H), for every $X'$ and $Y '$ in $\R^{2n'}$ and $Z''$ in
$\R^{2n''}$:
$$ K_{n' } (X'-Y') \left |  \pi ^{-n'}\int _{\R^{2n'}}
 F(Z', Z'')\   \Phi_{n' } (X' , Y' , Z')  d Z' \right  | \leq 9^{n'}
 N_ {n' } (F) , \leqno (3.11)$$
$$N_ {n' } (F) = \sum _{ (\alpha , \beta) \in I_2(n') }
 \Vert
 \partial _x^{\alpha} \partial _{\xi }^{\beta} F  \Vert
 _{L^{\infty}(\R^{2n})}$$
where $I_2(n')$ is the set of multi-indices $\alpha $
 such that $\alpha _j \leq 2$ and $\beta _j \leq 2$ for all $j\leq
 n$, and $\alpha_j = \beta_j =0$ for $j>n'$.

}

\bigskip

{\it Proof.} Let $I ( X',  Y', Z'')$ be the left side of (3.11).
Integrations by parts show that
for all  $X$ and $Y$ in  $\R^{2n'}$:
 $$ I ( X',  Y', Z'')  \leq \pi ^{-n'} \int _{\R^{2n'} }
\left | (L  F ) \left ( Z' + { X' + Y'\over 2}, Z'' \right ) \right
| e^{- | Z' |^2} dZ' $$
where $L$ is the differential operator defined by
$$ L = \prod _{j\leq n'} L_{z_j}  L_{\zeta_j}  \hskip 2cm
 L_{z_j}= \sum_{k=0}^2 a_k  ( z_j  )
\partial _{z_j}^k $$
$$ a_0 (z) = 3 - 4 z^2 \hskip 2cm a_1(z) = 4z \hskip 2cm a_2(z) = -1 .
$$
We get as a consequence that
$$ I ( X',  Y', Z'') \leq \sum _{(\alpha , \beta)\in I_2(n')}
\left \Vert \partial _x^{\alpha}
\partial _{\xi}^{\beta}   F \right \Vert _{L^{\infty}(\R^{2n})} \prod _{j\leq n'}
C_{\alpha_j} C_{\beta_j} $$
with
$$ C_k = \pi ^{-1/2} \int_{\R} |a_k(z)|e^{-z^2} dz \hskip 2cm k=0, 1 , 2 .$$
The formula (3.11) then follows from the fact that
 $\max (C_0 , C_1 , C_2 ) \leq 3$. \hfill $\carre$

\bigskip

 {\it End of the proof of Lemma  2.1.}
The subset  $E$ may be any subset of $\{ 1 , ... , n\}$, but we can
assume in the proof that $E = \{ 1 , ... , n'\}$ with $0 \leq n'
\leq n$. In this case, we use the above notations, and we set $x =
(x' , x'')$ with $x' = (x_1 , ... x_{n'})$ and $x'' = ( x_{n'+1} ,
... x_n)$, etc.

Using (3.6) and Lemma 3.3 to obtain an upper bound on the right side, one gets :
$$\Big | < Op ^{Weyl} ( e^{{1\over 4}\Delta ''} F) f, g> \Big |
\leq ... $$
$$ ... \leq 9^{n'} (2\pi )^{-2n'-n''} N_{n'} (F) \int K_{n'} (X' - Y')^{-1}
| < f , \Psi _{X', T''}> < \Psi _{Y', T'' }, g> | dX'dY' dT''$$
According to Schur's Lemma, this is smaller than
$$... \leq  9^{n'} (2\pi )^{-2n'-n''} N_{n'} (F) \Vert K_{n'} ^{-1}\Vert
_{L^1(\R^{2n'})}... $$
$$ ... \left [ \int |< f , \Psi _{X', T''}> |^2 dX' dT'' \right ]^{1/2} \
\left [ \int |<  \Psi _{Y', T''}, g> |^2 dY' dT'' \right ]^{1/2}.$$
Using (2.2), one shows that $ | < Op ^{Weyl} ( e^{{1\over 4}\Delta ''} F) f, g>|$ is smaller than
$$ \leq 9^{n'} (2\pi )^{-n'}  N_{n'} (F) \Vert K_{n'} ^{-1}\Vert
_{L^1(\R^{2n'})} \Vert f \Vert_{L^2(\R^n)}  \ \ \Vert g
\Vert_{L^2(\R^n)} .$$
Since $\Vert K_{n'} ^{-1}\Vert _{L^1(\R^{2n'})}  = \pi^{2n'}$, the former inequalities imply that
$$\Big | < Op ^{Weyl} ( e^{{1\over 4}\Delta ''} F) f, g>
\Big |  \leq \left ( {9 \pi \over 2} \right )^{n'}N_{n'} (F) \Vert f
\Vert_{L^2(\R^n)}  \ \ \Vert g \Vert_{L^2(\R^n)}. $$
 Lemma 2.1 holds if  $E = \{ 1 , ... , n'\}$,
a case which we are brought back to by a suitable permutation.

\bigskip

{\bf 4. Proof of Lemma 2.2}

\bigskip

Let $\Delta _j$ be the operator defined by   (2.6) and
$$ A_j = I - e^{{1\over 4}\Delta_j} \leqno (4.1)$$

\bigskip

{\bf Lemma 4.1.} {\it One can write
$$ A_j = B_j \partial _{ x_j} + C_j \partial _{ \xi_j}
 = D_j \Delta_j \leqno (4.2)$$
where the operators $B_j$, $C_j$ and $D_j$ are bounded in the space  $C_b$
of continuous bounded functions on  $\R^{2n}$. More precisely,

 $$ \Vert A_j \Vert _{{\cal L}(C_b)} \leq 2 \hskip 2cm
 \Vert B_j \Vert _{{\cal L}(C_b)} \leq  \pi ^{-1/2} \leqno (4.3)  $$
$$ \Vert C_j \Vert _{{\cal L}(C_b)} \leq  \pi ^{-1/2}
 \hskip 2cm
 \Vert D_j \Vert _{{\cal L}(C_b)} \leq 1/4.$$
}
\bigskip

{\it Proof.} The first inequality in (4.3) is standard.
The expression of the heat operator allows us to write
the first equality in
 (4.2) with
$$( B_j \varphi ) (x , \xi) = - \pi ^{-1} \int_{\R^2\times [0, 1]}
e^{- (u^2 + v^2)}u \varphi ( x + \theta u e_j , \xi + \theta v e_j
)dudv d\theta \leqno (4.4)$$
$$  ( C_j \varphi ) (x , \xi) = - \pi ^{-1} \int_{\R^2\times [0,
1]} e^{- (u^2 + v^2)}v \varphi ( x + \theta u e_j , \xi + \theta v
e_j )dudv d\theta \leqno (4.5)
$$
We deduce the bounds on the norms of
 $C_j$ and $D_j$ in
(4.3) from these inequalities. The last inequality  (4.2) and the bound on $D_j$
in (4.3) follow by integrating by parts in
  (4.4) and (4.5).  \hfill $\carre$

\bigskip

{\it End of the proof of Lemma 2.2.} For every multi-index
 $(\alpha
, \beta )$ in $I_2(E)$, one can rewrite the operator

$\partial_x^{\alpha}\partial_{\xi}^{\beta} T(E)$ as follows
$$\partial_x^{\alpha}\partial_{\xi}^{\beta} T(E) = \prod _{j\in E}
U_j \ \partial_x^{\alpha_j}\partial_{\xi}^{\beta_j}$$
with
$$ U_j =  \left \{ \matrix {A_j &{\rm if}& \alpha _j + \beta _j \geq 2\cr
B_j \partial _{x_j} + C_j \partial _{\xi _j} &{\rm if}& \alpha _j +
\beta _j = 1\cr D_j \Delta _j &{\rm if}& \alpha _j + \beta _j = 0\cr
} \right . $$

According to the bounds on the norms of
the operators $A_j$ to  $D_j$ given by Lemma 4.1, if  $F$
satisfies hypothesis (H) with $\rho_j = \delta_j \leq 1$, one has:
$$  \Vert \partial_x^{\alpha}\partial_{\xi}^{\beta} T(E) F
\Vert _{L^{\infty}(\R^{2n})} \leq M 2^{|E|}\prod _{j\in E} \rho_j^2
$$
 Since   $I_2(E)$ contains  exactly
$9^{|E|}$ elements,  this achieves the proof of Lemma 2.2.
 \hfill $\carre$

\bigskip

{\bf References.}

\smallskip

[A-J-N] L. Amour, L. Jager, J. Nourrigat, {\it Bounded Weyl
pseudodifferential operators in Fock space,} preprint,
arXiv:1209.2852, sept. 2012.

\smallskip

[C-V] A.P. Calder\'on, R. Vaillancourt, {\it A class of bounded
pseudo-differential operators,} Proc. Nat. Acad. Sci. U.S.A, {\bf
69}, (1972), 1185-1187.

\smallskip

[C-M] R. L. Coifman, Y. Meyer, {\it Au del\`a des op\'erateurs
pseudo-diff\'erentiels,} Ast\'erisque, {\bf 57}, 1978.

\smallskip

[C-R] M. Combescure, D. Robert, {\it Coherent states and
applications in mathematical physics,} Theoretical and Mathematical
Physics. Springer, Dordrecht, 2012.

\smallskip

[C] H. O. Cordes, {\it On compactness of commutators of
multiplications and convolutions, and boundedness of
pseudo-differential operators,} J. Funct. Anal, {\bf 18} (1975)
115-131.

\smallskip

[F] G. B. Folland, {\it Harmonic analysis in phase space.} Annals of
Mathematics Studies, {\bf 122}. Princeton University Press,
Princeton, NJ, 1989.

\smallskip

[HO] L. H\"ormander, {\it The analysis of linear partial
differential operators,} Volume III, Springer, 1985.

\smallskip

[HW] I. L. Hwang, {\it The $L^2$ boundedness of pseudo-differential
operators,} Trans. Amer. Math. Soc, {\bf 302} (1987) 55-76.

\smallskip

[LA1] B. Lascar, {\it Noyaux d'une classe d'op\'erateurs
pseudo-diff\'erentiels sur l'espace de Fock, et applications.}
S\'eminaire Paul Kr\'ee, 3e ann\'ee (1976-77), Equations aux
d\'eriv\'ees partielles en dimension infinie, Exp. No. 6, 43 pp.

\smallskip

[LA2] B. Lascar,  {\it Equations aux d\'eriv\'ees partielles en
dimension infinie.}  Vector space measures and applications (Proc.
Conf., Univ. Dublin, Dublin, 1977), I, pp. 286-313, Lecture Notes in
Math, {\bf 644}, Springer, Berlin, 1978.

\smallskip

[LA3] B. Lascar,  {\it Op\'erateurs pseudo-diff\'erentiels en
dimension infinie. Etude de l'hypoellipticit\'e de la
r\'esolubilit\'e dans des classes de fonctions h\"olderiennes et de
distributions pour des op\'erateurs pseudo-diff\'e\-rentiels
elliptiques,}  J. Analyse Math. {\bf 33} (1978), 39-104.

\smallskip

[LA4] B. Lascar,  {\it Une classe d'op\'erateurs elliptiques du
second ordre sur un espace de Hilbert,}  J. Funct. Anal. {\bf 35}
(1980), no. 3, 316-343.

\smallskip

[LA5] B.  Lascar, {\it  Op\'erateurs pseudo-diff\'erentiels en
dimension infinie. Applications. } C. R. Acad. Sci. Paris S\'er. A-B
{\bf 284} (1977), no. 13, A767-A769,

\smallskip

[LA6] B. Lascar, {\it  M\'ethodes $L^{2}$ pour des \'equations aux
d\'eriv\'ees partielles d\'ependant d'une infinit\'e de variables.}
 S\'eminaire Goulaouic-Schwartz (1975-1976) Equations aux
d\'eriv\'ees partielles et analyse fonctionnelle, Exp. No. 5, 11
pp. Centre Math., Ecole Polytech., Palaiseau, 1976.

\smallskip

[LA7] B.  Lascar, {\it  Propri\'et\'es locales d'espaces de type
Sobolev en dimension infinie. } Comm. Partial Differential Equations
{\bf 1} (1976), no. 6, 561-584.

\smallskip

[LA8] B. Lascar, {\it  Propri\'et\'es locales d'espaces de type
Sobolev en dimension infinie. } S\'eminaire Paul Kr\'ee, 1re ann\'ee
(1974-75), Equations aux d\'eriv\'ees partielles en dimension
infinie, Exp. No. 11, 16 pp.

\smallskip

[LA9] B. Lascar, {\it  Op\'erateurs pseudo-diff\'erentiels d'une
infinit\'e de variables, d'apr\`es M. I. Visik.} S\'eminaire Pierre
Lelong (Analyse), Ann\'ee 1973-1974, pp. 83–90. Lecture Notes in
Math., {\bf 474}, Springer, Berlin, 1975.

\smallskip

[LA10] B. Lascar, {\it  Propri\'et\'es d'espaces de Sobolev en
dimension infinie, } C. R. Acad. Sci. Paris S\'er. A-B {\bf 280}
(1975), no. 23, A1587-A1590.

\smallskip

[LER] N. Lerner, {\it Metrics on the phase space and non
self-adjoint pseudo-differential operators,} Birkh\"auser Springer,
2010.

\smallskip
[R] D. Robert, {\it Autour de l'approximation semi-classique},
Progress in Mathematics 68,
Birkh\"auser Boston, Inc., Boston, MA, 1987.

\smallskip

[U1] A. Unterberger, {\it  Oscillateur harmonique et op\'erateurs
pseudo-diff\'erentiels},   Ann. Inst. Fourier (Grenoble) {\bf 29}
(1979), no. 3, xi, 201-221.

\smallskip

[U2] A. Unterberger, {\it Les op\'erateurs m\'etadiff\'erentiels},
in Complex analysis, microlocal calculus and relativistic quantum
theory, Lecture Notes in Physics {\bf 126} (1980) 205-241.

\bigskip

\bigskip

Laboratoire de Math\'ematiques, FR CNRS 3399, EA 4535, Universit\'e de Reims
Champagne-Ardenne, Moulin de la Housse, B. P. 1039, F-51687
Reims, France,

{\it E-mail:} {\tt laurent.amour@univ-reims.fr }

{\it E-mail:} {\tt lisette.jager@univ-reims.fr}

{\it E-mail:} {\tt jean.nourrigat@univ-reims.fr}


 \end

 We have, for each multi-index $\alpha$ such that $\alpha_j \leq
2$ for all $j$ and for all $x\in \R^n$:
$$ |\partial^{\alpha} W(x)| \leq \prod _{j=1}^n ( C \lambda_j
)^{\alpha _j} \hskip 2cm C = 4 \max ( C_0, 1) $$
By the formula (2.9) of [C-S], we can write:
$$ \partial^{\alpha}e^{W(x)} = \alpha ! e^{W(x)} \sum _{(s , \beta , k)\in E(\alpha)}
 \prod _{j=1}^s {1 \over k_j ! } \left ( {\partial
^{\beta ^{(j)}}W(x) \over \beta^{(j)} ! } \right )^{k_j}$$
where $E( \alpha)$ is the set of elements $(s , \{ \beta ^{(1)} ,
... , \beta ^{(s)}, k_1 , ... , k_s \} )$ such that $s$ is an
integer between $1$ and $|\alpha|$, $\beta ^{(1)} , ... , \beta
^{(s)}$ are non zero multi-indices, all of them being different, and
$k_1 , ... , k_s$ are integers $\geq 1$, such that:
$$  \sum _{j= 1}^s k_j \beta^{(j)}= \alpha \leqno (1.8)$$
In the definition of $E( \alpha)$, we consider $\{ \beta ^{(1)} ,
... , \beta ^{(s)}, k_1 , ... , k_s \}$ as an unordered set, not as
a sequence. It follows that:
$$ |\partial^{\alpha}e^{W(x)}| \leq \alpha ! e^{W(x)}
|F(\alpha)|  \prod _{j=1}^n ( C \lambda_j )^{\alpha _j}   \leqno
(1.9)$$
where $F(\alpha)$ is the subset of  $E(\alpha)$ of elements $(s ,
\beta ^{(1)}  , ... , \beta ^{(s)}, k_1 , ... , k_s)$ such that none
of the $\partial ^{\beta^{(j)}}W$ vanishes identically. By the form
(1.7) of $W$, there is an integer $N$ with the following property:
for each point $j$, the number of multi-indices $\beta$ such that
$\beta _j \not= 0$, $\beta_k \leq 2$ for all $k$, and such that
$\partial ^{\beta}W$ does not vanish identically, is smaller than
$N$.  By (1.8), since $\alpha_j \leq 2$ for all $j$, all the numbers
$k_j$ in the definition of $E(\alpha)$ satisfy  $1 \leq k_j \leq 2$.
 Therefore, the number of elements of $F(\alpha)$  is smaller
than $( 4N)^{|\alpha|}$. Therefore, since $W \leq 0$ and $\alpha !
\leq 2^{|\alpha|} $, we have, with another $C$ depending only on
$C_0$:
-------------------------------------------------------------------------------------------
\smallskip

[C-S] G. M. Constantine, T. H. Savits, {\it A multivariate Faa di
Bruno formula with applications,} Trans. of the A.M.S, {\bf 348},
(2), (1996), 503-520.